\documentstyle{article}
\author{Johan Ernest Mebius\thanks{Delft University of Technology, Faculty of
Electrical Engineering, Mathematics and Computer Science, P.O.Box 5031, NL -- 2600 GA Delft,
The Netherlands, Phone +31.15.2783072, Fax +31.15.2786632}}
\title{A MATRIX--BASED PROOF OF THE\\QUATERNION REPRESENTATION THEOREM\\FOR FOUR-DIMENSIONAL ROTATIONS}
\date{September 2004}
\begin{document}
\maketitle

\begin{abstract}
\noindent
          To each 4x4 matrix of reals another 4x4 matrix is constructed,
          the so-called associate matrix. This associate matrix is shown
          to have rank 1 and norm 1 (considered as a 16D vector) if and
          only if the original matrix is a 4D rotation matrix.\\
          This rank-1 matrix is the dyadic product of a pair of 4D unit
          vectors, which are determined as a pair up to their signs.\\
          The leftmost factor (the column vector) consists of the components
          of the left quaternion and represents the left-isoclinic part
          of the 4D rotation. The rightmost factor (the row vector) likewise
          represents the right quaternion and the right-isoclinic part of the
          4D rotation.\\
          Finally the intrinsic geometrical meaning of this matrix-based proof
          is established by means of the usual similarity transformations.
\end{abstract}

\section*{Notations}

\begin{tabbing}
WWWWWWWWWWJ\=DUMMY      \kill
Bold--face symbols:     \>Geometrical entities\\
Normal--face symbols:   \>Their algebraic representations\\
                        \>\\
$A$:                    \>Van Elfrinkhof matrix (defined in Paragraph 4.2); Arbitrary 4D rotation matrix\\
$M$:                    \>Associate of matrix $A$ (defined in Paragraph 4.2)\\
$M_L$, $M_R$:           \>Matrices representing left-- and right--multiplication by a unit quaternion, respectively\\
                        \>(introduced in Paragraph 4.1)\\
$L = a + bi + cj + dk$: \>Unit quaternion appearing as a leftmost factor in quaternion multiplication\\
$P = u + xi + yj + zk$: \>Arbitrary 4D point represented as a quaternion\\
$R = p + qi + rj + sk$: \>Unit quaternion appearing as a rightmost factor in quaternion multiplication\\
$S$, $S_L$, $S_R$:      \>Similarity transformation matrix and its left-- and right--isoclinic components\\
                        \>(Similarity transformations restricted to rotations of coordinate system)\\
                        \>(introduced in Section 5)\\
\end{tabbing}

\section{Introduction}

The classical quaternion representation theorem for rotations in 4D Euclidean
space states that an arbitrary 4D rotation matrix is the product of a matrix
representing left--multiplication by a unit quaternion and a matrix representing
right--multiplication by a unit quaternion. This decomposition is unique up
to sign of the pair of component matrices.
\\
In this paper this theorem is proved as a result in the theory of matrices by
pure matrix means. As a matrix--theoretic theorem it has no {\em a priori}
geometrical meaning; for this one has to study the behaviour of its matrix
formulation under a predetermined class of similarity transformations.
\\
The proof in this paper is not the first one in the existing literature;
presumably {\sc Bouman} ([BOUM~1932]) was the first person to publish a
satisfactory proof of the representation theorem.
\\
\\
The application of quaternions to geometry is a well--known classical subject,
but in the opinion of the author it is also still a fertile research area.
Think only of the use of low--dimensional geometry and topology in certain
areas of theoretical physics. Or think of the hidden 4D rotational symmetry of the
{\sc Coulomb} field.
\\
The author's desire to make computer software for 3D and 4D geometry led to
the purely matrix--based proof of the quaternion representation theorem for 4D
rotations presented in this paper.

\section{Theorem and outline of proof} \label{Outline}

{\sc Theorem:} Each 4D rotation matrix can be decomposed in two
ways into a matrix representing left--multiplication by a unit quaternion
and a matrix representing right--multiplication by a unit quaternion.
These decompositions differ only in the signs of the component matrices.
\\
\\
{\sc Outline of proof:} Let $M_L$, $M_R$ be matrices representing
left-- and right--multiplication by a unit quaternion, respectively. Then
their product $A = M_LM_R$ is a 4D rotation matrix.
\\
Matrix $M_L$ is determined by four reals $a, b, c, d$ satisfying the relation
$a^2 + b^2 + c^2 + d^2 = 1$.
Likewise, matrix $M_R$ is determined by four reals $p, q, r, s$ satisfying the
relation $p^2 + q^2 + r^2 + s^2 = 1$.
\\
The 16 products $ap,\ aq,\ ar,\ as$, $\ldots$ , $dp,\ dq,\ dr,\ ds$ are arranged
into a matrix $M$, which has rank 1 and is easily expressed in the elements of
$A$. Let us in this paper denote it as the {\em associate matrix} of $A$.
\\
Conversely, given an arbitrary 4D rotation matrix $A$, one calculates its
associate matrix $M$ in the hope that it is a matrix of products $ap,\ aq,\ ar,\ as$,
$\ldots$ , $dp,\ dq,\ dr,\ ds$ which are not all zero. This hope is vindicated
by proving that $M$ has rank 1 whenever $A$ is a 4D rotation matrix.
\\
\\
The proof is completed by observing that the sum of the squares of the elements
of $M$ is unity, and concluding that two pairs of quadruples of reals
$a,\ b,\ c,\ d$; $p,\ q,\ r,\ s$ exist satisfying $a^2 + b^2 + c^2 + d^2 = 1$,
$p^2 + q^2 + r^2 + s^2 = 1$ and differing only in sign.

\section{Quaternions in 4D Euclidean geometry}

In this paper the 4D Euclidean space is provided with a Cartesian coordinate
system $OUXYZ$. Points are represented as column vectors $(u, x, y, z)^T$ (in
this paper denoted as the $R^4$ representation) or as quaternions
$u + xi + yj + zk$. In this connection the coordinate system is also denoted
as {\em O1IJK}.
\\
Both of these representations are based on an arbitrary choice of coordinate
system. To obtain a possible geometrical meaning of a matrix--algebraic
result, one has to prove that it is invariant under coordinate system
transformations ({\em similarity transformations}), which are in this paper
restricted to {\em rotations} because we are interested in Euclidean
properties, not in more general affine or projective properties.
\\
\\
In the sequel we have to do with geometrical objects as well as their
representations with respect to specific coordinate systems. For easy reading
object and representation are denoted by the same symbol; the object itself
in bold face, its algebraic representations in normal face, with decorations
added as needed.

\section{Proof of the representation theorem for 4D rotations}

Refer to Section \ref{Outline} above for the statement of the theorem and the
outline of its proof.

\subsection{Isoclinic 4D rotations}

Let ${\bf P}$ be an arbitrary 4D point, represented as a quaternion $P = u + xi + yj + zk$.
Let $L = a + bi + cj + dk$ and $R = p + qi + rj + sk$ be unit quaternions
($a^2 + b^2 + c^2 + d^2 = 1$, $p^2 + q^2 + r^2 + s^2 = 1$).
\\
Consider the left-- and right--multiplication mappings $M_L: P \to LP$ and
$M_R: P \to PR$. In the $R^4$ representation $M_L$ and $M_R$ are linear
mappings with matrices

\begin{equation}
\label{LandR}
M_L=\left [\begin {array}{rrrr} a&-b&-c&-d\\\noalign{\medskip}b&a&-d&c\\\noalign{\medskip}c&d&a&-b\\\noalign{\medskip}d&-
c&b&a\end {array}\right ],
\ M_R=\left [\begin {array}{rrrr} p&-q&-r&-s\\\noalign{\medskip}q&p&s&-r\\\noalign{\medskip}r&-s&p&q\\\noalign{\medskip}s&r
&-q&p\end {array}\right ].
\end{equation}

One easily proves that both $M_L$ and $M_R$ are orientation--preserving
isometries of $R^4$ with the origin of coordinates $O$ as a fixed point,
{\em i.\ e.} rotations of $R^4$ about $O$.
\\
\\
Both $M_L$ and $M_R$ have the property of rotating all half--lines originating
from $O$ through the same angle ($\arccos a$ for $M_L$ and $\arccos p$ for
$M_R$); such rotations are denoted as {\em isoclinic}.
There exists however a subtle difference between $M_L$ and $M_R$, which is
best illustrated by a specific example:
\\
Let $M_L$ and $M_R$ be left-- and right--multiplication by the quaternion
$Q_{\alpha} = \cos \alpha + i\sin \alpha$, respectively.
Then $M_L$ acts in both coordinate planes {\em 1I} and $JK$ as a rotation
through the angle $\alpha$, while $M_R$ acts in the {\em 1I} plane as a
rotation through $\alpha$ and in the $JK$ plane as a rotation through
$-\alpha$.
\\
The two kinds of isoclinic rotations are conveniently distinguished as
{\em left--} and {\em right--isoclinic}.
\\
\\
Conversely, an isoclinic 4D rotation about $O$ different from the non--rotation
$I$ and from the central reversion $-I$ is represented by either a left--multiplication
or a right--multiplication by a unique unit quaternion and so is either a
left-- or a right--isoclinic rotation. This theorem is not used in the sequel,
but it allows us to speak of "left-- or right--isoclinic" instead of the
cumbersome "represented by a unit quaternion left-- or right--multiplication".
\\
This theorem is presumably due to {\sc Robert S. Ball}; in [BALL~1889] the
author does not mention it explicitly as a theorem, but nevertheless gives a
proof. However, {\sc Ball}'s proof covers only the case in which the
off--diagonal elements of the given matrix are all nonzero. A complete proof
is given in [MEBI~1994].
\\
\\
In Section \ref{Intrinsic} we prove that the left-- and right--isocliny
properties and the rotation angle are independent of the choice of coordinate
system.

\subsection{General 4D rotations}

Let us study the composition of a left--isoclinic and a right--isoclinic
rotation. First apply $M_L$, then $M_R$ to a 4D point $P$: one obtains
$(LP)R$ in quaternion representation. First apply $M_R$, then $M_L$: one
obtains $L(PR)$. Quaternion multiplication is associative. Therefore left--
and right--isoclinic rotations are commutative, and we have
$M_LM_RP = M_RM_LP$ ($R^4$ representation) $ = LPR$ (quaternion representation).
In the $R^4$ representation the mapping $P \to LPR$ has the matrix

\begin{equation}
\label{LR}
M_LM_R =
\left [\begin {array}{rrrr} ap-bq-cr-ds&-aq-bp+cs-dr&-ar-bs-cp+dq&-as+br-cq-dp\\\noalign{\medskip}bp+aq-dr+cs&-bq+ap+
ds+cr&-br+as-dp-cq&-bs-ar-dq+cp\\\noalign{\medskip}cp+dq+ar-bs&-cq+dp-as-br&-cr+ds+ap+bq&-cs-dr+aq-bp
\\\noalign{\medskip}dp-cq+br+as&-dq-cp-bs+ar&-dr-cs+bp-aq&-ds+cr+bq+ap\end {array}\right ].
\end{equation}

{\sc Van Elfrinkhof} ([ELFR~1897]) was apparently the first person to treat
the relation between quaternion multiplications and 4D rotations in this
algebraic way. For this reason a matrix of the form of Eq.\ \ref{LR} with
\mbox{$a^2 + b^2 + c^2 + d^2 = 1$} and \mbox{$p^2 + q^2 + r^2 + s^2 = 1$}
is in this paper denoted as a {\em Van Elfrinkhof matrix}.
\\
\\
A Van Elfrinkhof matrix

\begin{equation}
\label{A}
A =
\left[
\begin{array}{rrrr}
\noalign{\medskip}  a_{00}  &   a_{01}  &   a_{02}  &   a_{03}  \\
\noalign{\medskip}  a_{10}  &   a_{11}  &   a_{12}  &   a_{13}  \\
\noalign{\medskip}  a_{20}  &   a_{21}  &   a_{22}  &   a_{23}  \\
\noalign{\medskip}  a_{30}  &   a_{31}  &   a_{32}  &   a_{33}
\end{array}
\right]
\end{equation}
\\
is readily interpreted as a set of 16 linear equations in the 16 unknowns
$ap,\ aq,\ ar,\ as$, $\ldots$ , $dp,\ dq,\ dr,\ ds$.
With slightly more work than just an inspection of the plus and minus signs
in Eq.\ \ref{LR} one obtains these unknowns, arranged as a
matrix, which is in turn written as a dyadic product:

\begin{equation}
\label{M1}
M =
\left[
\begin{array}{rrrr}
\noalign{\medskip}  ap  &   aq  &   ar  &   as  \\
\noalign{\medskip}  bp  &   bq  &   br  &   bs  \\
\noalign{\medskip}  cp  &   cq  &   cr  &   cs  \\
\noalign{\medskip}  dp  &   dq  &   dr  &   ds
\end{array}
\right]
=
\left[
\begin{array}{rrrr}
\noalign{\medskip}  a  \\
\noalign{\medskip}  b  \\
\noalign{\medskip}  c  \\
\noalign{\medskip}  d
\end{array}
\right]
\begin{array}{r}
\left[
\begin{array}{rrrr}
\noalign{\medskip}  p  &   q  &   r  &   s
\end{array}
\right]
\\
\begin{array}{rrrr}
\   & \   & \   & \   \\
\   & \   & \   & \   \\
\   & \   & \   & \   \\
\   & \   & \   & \
\end{array}
\end{array}
\end{equation}

\begin{equation}
\label{M2}
= \frac{1}{4}
\left[
\begin{array}{rrrr}
\noalign{\medskip} a_{00}+a_{11}+a_{22}+a_{33} & +a_{10}-a_{01}-a_{32}+a_{23} & +a_{20}+a_{31}-a_{02}-a_{13} & +a_{30}-a_{21}+a_{12}-a_{03} \\
\noalign{\medskip} a_{10}-a_{01}+a_{32}-a_{23} & -a_{00}-a_{11}+a_{22}+a_{33} & +a_{30}-a_{21}-a_{12}+a_{03} & -a_{20}-a_{31}-a_{02}-a_{13} \\
\noalign{\medskip} a_{20}-a_{31}-a_{02}+a_{13} & -a_{30}-a_{21}-a_{12}-a_{03} & -a_{00}+a_{11}-a_{22}+a_{33} & +a_{10}+a_{01}-a_{32}-a_{23} \\
\noalign{\medskip} a_{30}+a_{21}-a_{12}-a_{03} & +a_{20}-a_{31}+a_{02}-a_{13} & -a_{10}-a_{01}-a_{32}-a_{23} & -a_{00}+a_{11}+a_{22}-a_{33}
\end{array}
\right].
\end{equation}
\\
Matrix $M$ contains at least one nonzero element and consists of mutually
proportional columns and therefore has rank 1. The sum of the squares of its
elements is $(a^2 + b^2 + c^2 + d^2)(p^2 + q^2 + r^2 + s^2) = 1$.
In this paper $M$ is denoted as the {\em associate matrix} of the matrix $A$.

\subsection{The general rotation as a product of isoclinic rotations}

In this paragraph we prove that any 4D rotation matrix can be decomposed
into a left-- and a right--isoclinic rotation matrix in two ways, differing
only in sign. If we can find unit quaternions $L = a + bi + cj + dk$ and
$R = p + qi + rj + sk$ such that the given 4D rotation matrix can be written
as a Van Elfrinkhof matrix, then we are done.
\\
\\
Let $A$ (Eq.\ \ref{A}) be an arbitrary 4D rotation matrix, and let M, calculated
according to Eq.\ \ref{M2}, be its associate matrix.
\\
Easy and somewhat tedious calculations show that the sum of the squares of
the elements of $M$ is unity, and that all 2--by--2 minors of $M$
are zero. In these calculations one uses the orthogonality of $A$ and the
equality of complementary 2--by--2 minors of 4D rotation matrices.
\\
The equality of complementary 2--by--2 minors follows from the general
theorem stating that complementary minors of an orthogonal matrix $B$ of
arbitrary order are equal or opposite according to Det $B$ equalling $+1$
or $-1$. ([CESA~1904], Art.\ 71, p.\ 54)
\\
It follows that $M$ has rank 1, and there exist reals $a,\ b,\ c,\ d$; $p,\ q,\ r,\ s$
satisfying $a^2 + b^2 + c^2 + d^2 = 1$, $p^2 + q^2 + r^2 + s^2 = 1$ such that
$M$ actually obeys Eq.\ \ref{M1}
\\
This yields two pairs of unit quaternions $L = a + bi + cj + dk$; $R = p + qi + rj + sk$,
differing only in sign. Finally, Eq.\ \ref{LandR} gives the left-- and
right--isoclinic components, and Eq.\ \ref{LR} the original 4D rotation matrix
as a Van Elfrinkhof matrix.
\\
\\
{\em Remark:} It is sufficient to prove that the second, third and fourth
columns of $M$ are proportional to its first column. So only nine
2--by--2 minors need to be calculated, not all 36 of them. The Appendix
shows one of these calculations.

\section{Independence of choice of coordinate system} \label{Intrinsic}

In this section we prove that the properties of left-- and right--isocliny
and the rotation angle of an isoclinic matrix are invariant under arbitrary
rotations of coordinate system. This section is in part a refresher on
dispacement transformations, similarity transformations, and their interplay.
\\
\\
Let ${\bf A}$ be a 4D rotation with centre $O$. Let $A$ be its matrix
referred to a Cartesian coordinate system $OUXYZ$. Let $OU'X'Y'Z'$ a be second
Cartesian coordinate system, originating from the first one by a rotation
${\bf S}$ about $O$. Rotation ${\bf A}$ is a displacement transformation;
rotation ${\bf S}$ is a similarity transformation. Let $S$ be the matrix of
${\bf S}$ referred to $OUXYZ$ and let $S_L$, $S_R$ be its left-- and
right--isoclinic components. It does not matter which one of the two possible
decompositions is taken.
\\
Furthermore, let ${\bf p}$ be an arbitrary 4D point, $p = (u,\ x,\ y,\ z)^T$,
$p' = (u',\ x',\ y',\ z')^T$ its coordinates referred to $OUXYZ$,
$OU'X'Y'Z'$, respectively. Then $p' = S^{-1}p$ expresses the new coordinates
of ${\bf p}$ in the old ones, and $p = Sp'$ the old coordinates in the new
ones. Let ${\bf Q} = {\bf Ap}$  point ${\bf p}$ rotated by ${\bf A}$.
\\
How does similarity transformation ${\bf S}$ affect rotation matrix $A$? To
find this out, begin at new coordinates $p'$, transform to old coordinates
$p = Sp'$, apply rotation matrix $A$ to obtain $Q = Ap = ASp'$, and finally
transform back to new coordinates, ending at $Q' = A'p' = S^{-1}ASp'$.
We conclude that the matrix $A'$ of ${\bf A}$ referred to $OU'X'Y'Z'$ equals
$A' = S^{-1}AS$.
\\
\\
Now we come to prove the intrinsic geometrical nature of the isocliny
properties of 4D rotations. Consider the general rotational similarity
transformation of the general 4D rotation:

\begin{equation}
\label{simtrans}
A' = S^{-1}AS = S_R\! ^{-1} S_L\! ^{-1}. A_L A_R. S_L S_R = S_L\! ^{-1} A_L S_L. S_R\! ^{-1} A_R S_R,
\end{equation}
\\
where the commutativity of left-- and right--isoclinic matrices was applied
to sweep all left--isoclinic components to the one side while retaining their
order, and the right--isoclinic components to the other side while retaining
their order too.
\\
Applying Eq.\ \ref{simtrans} to a left--isoclinic matrix $A_L$ (set $A_R = I$)
one obtains the transformation formula
\[
A_L' = S^{-1}A_LS = S_L\! ^{-1} A_L S_L,
\]
which shows that $A_L'$, being the product of left--isoclinic factors,
is left--isoclinic. One may attribute the property of
left--isocliny to the rotation ${\bf A_L}$ represented by matrix $A_L$ and
all its transforms. This establishes the geometrical nature of this property.
In other words, one may speak of left--isoclinic rotations in their own right,
irrespective of particular matrix representations.
\\
Likewise, for a right--isoclinic matrix one obtains
\[
A_R' = S^{-1}A_RS = S_R\! ^{-1} A_R S_R,
\]
which shows that $A_R'$ is right--isoclinic. Here too the property of right--
isocliny may be attributed to the rotation ${\bf A_R}$ represented by $A_R$
and all its transforms.
\\
In any rotation matrix of any order the main diagonal elements are the cosines
of the angles through which the coordinate axes are rotated. Therefore, in an
isoclinic matrix of either kind the main diagonal elements are equal.
Their sum, the trace of the matrix, is invariant under similarity transformations.
It follows that the rotation angle of an isoclinic matrix is also an invariant.
So we may speak of the rotation angle of an isoclinic rotation, irrespective
of particular matrix representations.

\section*{References}

[BALL~1889] {\sc Robert S. Ball}: {\em Theoretische Mechanik starrer Systeme}.
Herausgeber: {\sc Harry Gravelius}. Berlin: Georg Reimer, 1889
\\

[BOUM~1932] {\sc J.\ Bouman}: Over quaternionen en hunne toepassing in de
meetkunde der vierdimensionale ruimte. {\em Nieuw Archief voor Wiskunde,
reeks 2, deel XVII, stukken 3 en 4}, 1932, p. 240-266
\\

[CESA~1904] {\sc E.\ Ces\`aro}: {\em Elementares Lehrbuch der algebraischen Analysis
und der Infinitesimalrechnung}. Leipzig: B.\ G.\ Teubner, 1904
\\

[ELFR~1897] {\sc L.\ van Elfrinkhof}: Eene eigenschap van de orthogonale
substitutie van de vierde orde. {\em Handelingen van het 6$^e$Nederlandsch
Natuurkundig en Geneeskundig Congres, Delft, 1897}
\\

[MEBI~1994] {\sc Johan Ernest Mebius}: {\em Applications of quaternions to
dynamical simulation, computer graphics and biomechanics}. Ph.D. Thesis
Delft University of Technology, Delft, 1994

\section*{Acknowledgement}

The author wishes to thank {\sc Th.\ H.\ M.\ Smits} for his useful comments
and suggestions.

\section*{Appendix}

The sum of the squares of the elements of $M$ in \mbox{Eq.\ \ref{M2}} is found to be
\[
4(a_{00}^{2}+a_{01}^{2}+a_{02}^{2}+a_{03}^{2}
+a_{10}^{2}+a_{11}^{2}+a_{12}^{2}+a_{13}^{2}
+a_{20}^{2}+a_{21}^{2}+a_{22}^{2}+a_{23}^{2}
+a_{30}^{2}+a_{31}^{2}+a_{32}^{2}+a_{33}^{2})/16,
\]
which equals unity because of the orthogonality of $A$.
\\
\\
All 2--by--2 minors of $M$ are readily shown to be zero. For instance,
\[
16(m_{00}m_{11} - m_{10}m_{01}) = (a_{00}+a_{11}+a_{22}+a_{33})(-a_{00}-a_{11}+a_{22}+a_{33}) - (a_{10}-a_{01}+a_{32}-a_{23})(a_{10}-a_{01}-a_{32}+a_{23})
\]
\[
= (a_{22} + a_{33})^2 - (a_{00} + a_{11})^2 - (a_{10} - a_{01})^2 + (a_{32} - a_{23})^2
\]
\begin{equation}
\label{Minor}
= a_{22}^2 + a_{33}^2 + a_{32}^2 + a_{23}^2 - a_{00}^2 - a_{11}^2 - a_{10}^2 - a_{01}^2
+ 2 (a_{22}a_{33} - a_{32}a_{23} - a_{00}a_{11} + a_{10}a_{01});
\end{equation}
the bilinear terms cancel because in the rotation matrix $A$ the
complementary minors $a_{22}a_{33} - a_{32}a_{23}$ and $a_{00}a_{11} - a_{10}a_{01}$ are equal.
\\
\\
As regards the quadratic terms, we have
\[
(a_{22}^2 + a_{33}^2 + a_{32}^2 + a_{23}^2) + (a_{02}^2 + a_{12}^2 + a_{03}^2 + a_{13}^2) = 2,
\]
\[
(a_{00}^2 + a_{11}^2 + a_{10}^2 + a_{01}^2) + (a_{02}^2 + a_{12}^2 + a_{03}^2 + a_{13}^2) = 2,
\]
therefore the quadratic terms in Eq.\ \ref{Minor} cancel too.

\end{document}